\documentclass[12pt]{article}
\usepackage{mathrsfs}
\usepackage{amssymb}
\usepackage{amsmath}
\usepackage{graphics}
\usepackage{epsfig}
\newtheorem{Corollary}{Corollary}[part]
\newtheorem{Definition}{Definition}[part]

\newtheorem{Lemma}{Lemma}[part]
\newtheorem{Proposition}{Proposition}[part]
\newtheorem{Remark}{Remark}[part]
\newtheorem{Theorem}{Theorem}[part]

\numberwithin{Assumption}{section}
\numberwithin{Corollary}{section}
\numberwithin{Definition}{section}
\numberwithin{equation}{section} \numberwithin{Example}{section}
\numberwithin{Lemma}{section} \numberwithin{Proposition}{section}
\numberwithin{Remark}{section} \numberwithin{Theorem}{section}

\def \tsm{\mathop{\textstyle\sum}\limits}

\def \wt{\widetilde}

\def \R{{\mathbb{R}}}
\def \E{{\mathbb{E}}}
\def \P{{\mathbb{P}}}
\def \Q{{\mathbb{Q}}}

\def \Dc{\mathcal{D}}

\def \Fc{\mathcal{F}}

\def \Lc{\mathscr{L}}

\def \Mc{\mathcal{M}}
\def \Mcc{\mathscr{M}}

\def \id{{\mathbf{1}}}
\def \qed{{\hbox{ }\hfill$\Box$}}
\def \tr{\mathrm{tr}}

\def\xa{\alpha}

\def\xd{\delta}
\def\xD{\Delta}

\def\xs{\sigma}
\def\xT{\Theta}

\def\xvt{\vartheta}

\def\xO{\Omega}

\def\reff#1{{\rm(\ref{#1})}}

\def\var{\mathrm{Var}}

\begin{document}

\title{\large\bf Mean-variance Hedging in the Discontinuous Case}
\author{{\large\bf Jianming Xia}
\footnote{Supported by the National Natural Science Foundation of
China under grant 10571167. Email: xia@amss.ac.cn.}\\
\small Academy of Mathematics and Systems Science, Chinese Academy
of Sciences\\
\small P.O. Box 2734, Beijing 100080, P.R. China }
\date{}
\maketitle

\begin{abstract}

The results on the mean-variance hedging problem in Gouri\'eroux,
Laurent and Pham (1998), Rheinl\"ander and Schweizer (1997) and Arai
(2005) are extended to discontinuous semimartingale models. When the
num\'eraire method is used, we only assume the Radon-Nikodym
derivative of the variance-optimal signed martingale measure (VSMM)
is non-zero almost surely (but may be strictly negative). When
discussing the relation between the solutions and the
Galtchouk-Kunita-Watanabe decompositions under the VSMM, we only
assume the VSMM is equivalent to the reference probability.

\

{\bf JEL Classification:} G10

{\bf 2000 Mathematics Subject Classification:} 91B28, 60H05, 60G48

{\bf Key words and phrases:} Mean-variance hedging, variance-optimal
martingale measure, num\'eraire, Galtchouk-Kunita-Watanabe
decomposition

\end{abstract}

\newpage

\section{Introduction}\label{sec-intro}
\setcounter{equation}{0} \setcounter{Assumption}{0}
\setcounter{Theorem}{0} \setcounter{Proposition}{0}
\setcounter{Corollary}{0} \setcounter{Lemma}{0}
\setcounter{Definition}{0} \setcounter{Remark}{0}

Let $S$ be a semimartingale and $\xT$ a family of some
$S$-integrable predictable processes $\xvt$ such that the stochastic
integral $G_T(\xvt):=\int_0^T\xvt_t\,dS_t\in L^2(\P)$, where $T$ is
a positive time horizon, then $G_T(\xT):=\{G_T(\xvt): \xvt\in\xT\}$
is a subspace of $L^2(\P)$. The problem of mean-variance hedging is
to approximate any contingent claim, i.e., any random variable $H\in
L^2(\P)$ by the elements in $G_T(\xT)$. In order to guarantee the
existence of the solution of such a problem, the working space $\xT$
of admissible strategies should be chosen such that $G_T(\xT)$ is
closed in $L^2(\P)$.

In the existing literature, the space $\xT$ usually consists of all
$S$-integrable predictable processes $\xvt$ such that the stochastic
integral $G(\xvt):=\int\xvt\, dS$ is a square-integrable
semimartingale. If $S$ is a (local) martingale, the closedness holds
true by the definition of stochastic integration. If $S$ is only a
semimartingale, additional assumptions must be imposed to ensure the
closedness. For a continuous semimartinagle, Delbaen et al. (1997)
established necessary and sufficient conditions for the closedness.
For further results along this line, see Grandits and Krawczyk
(1998) and Choulli et al. (1998, 1999). When the problem of
mean-variance hedging is studied, $G_T(\xT)$ is usually assumed to
be closed in $L^2(\P)$ explicitly or implicitly under additional
conditions, see Schweizer (1996), Rheinl\"ander and Schweizer (1997)
(RS 1997, for short), Hou and Karatzas (2004) and Arai (2005), among
others. But all these additional conditions imposed on $S$ are
rather strong.

On the other hand, Delbaen and Schachermayer (1996b) defined the
working space starting from ``simple" strategies and discussed the
duality relation between attainable claims (by admissible
strategies) and equivalent martingale measures. The space chosen by
them automatically has the $L^2(\P)$-closedness. Inspired by this,
for continuous semimartingale models, Gouri\'eroux, Laurent and Pham
(1998) (GLP 1998, for short) dealt with the mean-variance hedging
problem, using the same working space as in Delbaen and
Schachermayer (1996b).

It is well known that the variance-optimal signed martingale measure
(VSMM, for short) plays an important role in studying the
mean-variance hedging problem. For a continuous semimartingale
model, it turns out that VSMM is equivalent to the reference
probability measure, see Delbaen and Schachermayer (1996a). But in
general, the VSMM is only a signed measure and therefore the set of
equivalent martingale measures is not enough. Inspired by this and
by a similar way of Delbaen and Schachermayer (1996b), when
discussing Markowitz's portfolio selection problem, Xia and Yan
(2006) defined a space of admissible strategies which has the
$L^2(\P)$-closedness and has the duality relation to signed
martingale measures (rather than only to equivalent ones). In an
independent work, \u{C}ern\'y and Kallsen (2005) also observed this
fact and chose the same space.

The aim of this paper is to investigate the mean-variance hedging
problem for discontinuous models within the working space of Xia and
Yan (2006). The results of RS 1997, GLP 1998 and Arai (2005) are
extent to our settings. RS 1997 and GLP 1998 dealt with the
continuous semimartingale model. Using a change of num\'eraire and a
change of measure, GLP 1998 reduced the problem to a martingale
framework. RS 1997 used the Galtchouk-Kunita-Watanabe decomposition
(GKW decomposition, for short) under the VSMM and obtained a
solution of feedback form. They also discussed the relation between
their solutions and those of GLP 1998. Arai (2005) extended the
results of GLP 1998 and RS 1997 to the discontinuous case under some
additional assumptions on the VSMM: the VSMM is equivalent to the
reference probability, the corresponding density process $Z$
satisfies the reverse H\"older inequality and another condition on
the jump of $Z$. But in our paper, when the num\'eraire method is
used, we only assume the Radon-Nikodym derivative of the VSMM is
non-zero almost surely (but may be strictly negative). When
discussing the relation between the solutions and the GKW
decompositions under the VSMM, we only assume the VSMM is equivalent
to the reference probability. The other rather restrictive
conditions such as the reverse H\"older inequality are removed here.

\section{The market model}\label{sec-set}
\setcounter{equation}{0} \setcounter{Assumption}{0}
\setcounter{Theorem}{0} \setcounter{Proposition}{0}
\setcounter{Corollary}{0} \setcounter{Lemma}{0}
\setcounter{Definition}{0} \setcounter{Remark}{0}

Let $(\Omega, \Fc, (\Fc_t)_{0\le t\le T}, \P)$ be a filtered
probability space satisfying the usual conditions, where
$\Fc_0=\sigma\{\emptyset, \Omega\}$, $\Fc_T=\Fc$, and $T$ is a
positive time horizon. Throughout this paper, $L^2(\xO,\Fc,\P)$ is
abbreviated as $L^2(\P)$. For any $a,b\in\R$, we denote $a\vee
b=\max\{a,b\}$ and $a\land b=\min\{a,b\}$. All vectors are column
vectors and the transposition of a vector is denoted by $x^\tr$. For
any $x,y\in\R^d$, the inner product of $x$ and $y$ is $x^\tr y$ and
the Euclidean norm of $x$ is $|x|:=\sqrt{x^\tr x}$.


\subsection{Simple strategies and signed martingale measures}

In this subsection, we first introduce the definitions of simple
strategies and signed martingale measures and then present some
existing results.

\begin{Definition}\label{def-local-l2}
The family $\Lc^2(\P)$ consists of all $\R^d$-valued
$(\Fc_t)$-progressively  measurable processes $S$ such that $\{S_U:
U \mbox{ stopping time}\}\subset L^2(\P)$. The family
$\Lc^2_{\mathrm{loc}}(\P)$ consists of all $\R^d$-valued
$(\Fc_t)$-progressively  measurable processes $S$ such that there
exists a sequence $(U_n)_{n\ge1}$ of localizing stopping times
increasing to $T$ such that, for each $n\ge1$, the stopped process
$S^{U_n}\in\Lc^2(\P)$.
\end{Definition}

We make the following standing hypothesis on an $\R^d$-value process
$S$, which models the (discounted) price processes of the risky
assets:

\begin{description}
    \item[(H0)] $S\in\Lc^2_{\mathrm{loc}}(\P)$.
\end{description}

\begin{Remark}{\rm Under (H0), $S$ is not necessarily a semimartingale.
We in fact don't even assume $S$ is optional at the moment.}
\end{Remark}

\begin{Definition} We say a process $\xvt$ is a simple trading strategy if $\xvt$
has a form
\begin{eqnarray}\label{trade-simple}\xvt=\sum_{i=1}^nh_i\id_{]\!]T_{1i},T_{2i}]\!]},\end{eqnarray}
where, for each $i=1,\dots,n$, $T_{1i}\le T_{2i}$ are stopping times
such that $S^{T_{2i}}\in\Lc^2_{\mathrm{loc}}(\P)$ and $h_i$ is
bounded $\R^d$-valued $\Fc_{T_{1i}}$-measurable. The space $\xT^s$
consists of all simple trading strategies. For any $\xvt\in\Theta^s$
having form \reff{trade-simple}, the stochastic integral of $\xvt$
with respect to $S$ is $$G_t(\xvt):=(\xvt\bullet S)_t
=\sum_{i=1}^nh^{\mathrm{tr}}_i(S_{T_{2i}\land t}-S_{T_{1i}\land
t}).$$ Obviously, $G_T(\Theta^s):=\{G_T(\xvt): \xvt\in\Theta^s\}$ is
a subspace in $L^2(\P)$.
\end{Definition}

We will use the following notations:
\begin{itemize}
  \item $\Dc^s:=\{g\in L^2(\P): \E[gf]=0 \mbox{ for
all }f\in G_T(\Theta^s) \mbox{ and } \E[g]=1\}$;
  \item $\Dc^e:=\{g\in\Dc^s: g>0 \mbox{ a.s.}\}$;
  \item For any $g\in\Dc^s$, $\Q^g$ is the signed measure on $(\xO,\Fc)$ defined by ${d\Q^g\over d\P}=g$;
  \item For any $g\in\Dc^s$, $Z^g$ is the RCLL version of the martingale $(\E[g|\Fc_t])$;
  \item $\Mc^s:=\{\Q^g: g\in\Dc^s\}$ and $\Mc^e=\{\Q^g:
  g\in\Dc^e\}$.
\end{itemize}
It is clear that $\Dc^s$ is convex and closed in $L^2(\P)$. For each
$g\in\Dc^s$ and each $i=1,\dots,d$, $S^iZ^g$ is a local martingale.


\begin{Definition} Any element in $\Mc^s$ (resp. $\Mc^e$) is called a signed
(resp. equivalent) martingale measure for $S$. \end{Definition}

Throughout this paper, the closure $\overline{\{\cdots\}}$ refers to
the $L^2(\P)$-norm. Then we have the following easy lemma, see,
e.g., Lemma 2.1 of Delbaen and Schachermayer (1996a).

\begin{Lemma}\label{lma-ds-1}
Under assumption (H0), we have:
\begin{description}
\item[(a)] $\Mc^s\ne\emptyset\Longleftrightarrow
    1\notin \overline{G_T(\Theta^s)}$;
\item[(b)] For any $g\in L^2(\P)$,
    \begin{eqnarray*}g\in\Dc^s
    &\Longleftrightarrow&\E[g]=1\mbox{ and  }\E[gf]=0\mbox{ for all }f\in
    \overline{G_T(\Theta^s)}.\end{eqnarray*}
\end{description}
\end{Lemma}

Lemma \ref{lma-ds-1}(b) and the following lemma, which goes back to
Lemma 2.2 of Xia and Yan (2006), give the bipolar relation between
$\overline{G_T(\Theta^s)}$ and $\Dc^s$.

\begin{Lemma}\label{lma-bipolar} Assume (H0) and $\Mc^s\ne\emptyset$,
then we have:
\begin{eqnarray}\label{eq-bipolar}f\in\overline{G_T(\Theta^s)}
\Longleftrightarrow f\in L^2(\P) \mbox{ and } \E[fg]=0 \mbox{ for
all }g\in\Dc^s.\label{polar}\end{eqnarray}
\end{Lemma}

\subsection{Admissible trading strategies}

Subsequently, we always assume $S$ satisfies the following
condition:

\begin{description}
  \item[(H1)] $S$ is an $\R^d$-valued RCLL semimartingale and
  $S\in\Lc^2_{\mathrm{loc}}(\P)$.
\end{description}

The stochastic integral of a predictable process $\vartheta$ with
respect to a semimartingale $X$ is denoted as $\int \vartheta\,dX$
or $\vartheta\bullet X$. We denote by $\Lc(X)$ the set of all
$X$-integrable predictable processes. For the theory of stochastic
integration we refer to Jacod (1979), and Jacod and Shiryaev (1987);
particularly, for vector stochastic integrals, see Jacod (1980), and
Shiryaev and Cherny (2002). Following Xia and Yan (2006) (see also
an independent work of \v{C}ern\'y and Kallsen (2005)), we give the
definition of admissible strategies below.

\begin{Definition}\label{def-adm}
An admissible trading strategy is a process $\vartheta\in\Lc(S)$
such that $G_T(\xvt):=(\xvt\bullet S)_T\in\overline{G_T(\Theta^s)}$.
The space $\Theta$ consists of all admissible trading strategies and
$G_T(\xT):=\{G_T(\xvt):\ \xvt\in\xT\}$.

\end{Definition}\begin{Definition} For any $\xvt\in\Lc(S)$, $(\xvt^j)$ is a sequence
of simple strategies approximating to $\xvt$, if
$(\xvt^j)\subset\Theta^s$ and $G_T(\xvt^j)\to G_T(\xvt)$ in
$L^2(\P)$.\end{Definition}

By definitions, for any $\xvt\in\Lc(S)$, $\vartheta$ is admissible
if and only if it allows an approximating sequence of simple
strategies. Lemma \ref{lma-bipolar} yields

\begin{Lemma}\label{lma-adm} Under assumptions (H1) and that $\Mc^s\ne\emptyset$,
for any $\xvt\in\Lc(S)$, $\xvt\in\Theta$ if and only if $\xvt$
satisfies the following condition:
\begin{eqnarray}\left\{%
\begin{array}{ll}
    G_T(\xvt):=\int_0^T\xvt_t\,dS_t\in L^2(\P)\\
    \E[G_T(\xvt)g]=0\quad \mbox{for all }g\in\Dc^s.
\end{array}%
\right.\label{cond-admissible}
\end{eqnarray}
\end{Lemma}

Obviously, $\Theta^s\subset \Theta$, $G_T(\Theta^s)\subset G_T(\xT)$
and $G_T(\xT)\subset \overline{G_T(\Theta^s)}$. Thus we have
$\overline{G_T(\Theta^s)}=\overline{G_T(\xT)}$. Furthermore, the
following theorem, which goes back to Xia and Yan (2006, Theorem 2.1
and Remark 2.3), shows that $G_T(\xT)$ is automatically closed in
$L^2(\P)$, if we assume in addition that
\begin{description}
    \item[(H2)] $\Mc^e\ne\emptyset$.
\end{description}

\begin{Theorem}\label{thm-adm-sim}
Under assumptions (H1) and (H2), we have:
\begin{description}
    \item[(a)] For any $f\in\overline{G_T(\Theta^s)}$, there exists a
$\xvt\in\xT$ such that $f=\int_0^T\xvt_t\, dS_t$ and $\int\xvt\, dS$
is a uniformly integrable $\Q$-martingale for each $\Q\in\Mc^e$;
    \item[(b)] $\overline{G_T(\Theta^s)}=G_T(\xT)$.
\end{description}
\end{Theorem}

\begin{Remark}\label{rmk-vtq}{\rm
Under assumptions of the previous theorem, if $\xvt\in\Theta$
satisfies the conditions in (a) and $(\xvt^j)\subset\Theta^s$ is an
approximating sequence, then for any $\Q\in\Mc^e$, $(\xvt^j\bullet
S)_T\to(\xvt\bullet S)_T$ in $L^1(\Q)$. On the other hand, for any
$\vartheta^j\in\Theta^s$, $\int\xvt^j\, dS$ is a uniformly
integrable $\Q$-martingale for each $\Q\in\Mc^e$ (see Lemma
\ref{lma-xth-s} below). Thus for any $\Q\in\Mc^e$ and any stopping
time $\tau$,
$$(\xvt^j\bullet S)_{\tau}=\E_\Q[(\xvt^j\bullet
S)_T|\Fc_{\tau}]\overset{L^1(\Q)}{-\!-\!\!\!\longrightarrow}\E_\Q[(\xvt\bullet
S)_T|\Fc_{\tau}]=(\xvt\bullet S)_{\tau},$$ which implies
\begin{eqnarray}\label{lim-vt-2}(\xvt^j\bullet
S)_{\tau}\overset{\P}{\longrightarrow}(\xvt\bullet
S)_{\tau}\quad\mbox{for any stopping time }\tau.\end{eqnarray}
This fact will be used in proving Theorem
\ref{thm-theta-u}.}\end{Remark}

\begin{Remark}{\rm Kreps-Yan theorem (see, e.g., Schachermayer 2005) yields
$$\Mc^e\ne\emptyset
    \Longleftrightarrow\overline{G_T(\Theta^s)-L^2_+(\P)}\bigcap L^2_+(\P)=\{0\}$$
}\end{Remark}

\begin{Definition} The space $\Theta^u$ consists of all processes $\vartheta\in\Lc(S)$ satisfying
\begin{eqnarray}
\left\{%
\begin{array}{ll}
    G_T(\xvt):=\int_0^T\xvt_t\,dS_t\in L^2(\P)\\
    (\xvt\bullet S)Z^g \mbox{ is a uniformly integrable martingale for each
    }g\in\Dc^s.
\end{array}%
\right.\label{cond-admissible-s}
\end{eqnarray}
\end{Definition}

\begin{Theorem}\label{thm-theta-u}
Under assumptions (H1) and (H2), $G_T(\Theta)=G_T(\Theta^u)$.
\end{Theorem}

{\bf Proof.}  ``$\supset$" is clear. Conversely, for any $f\in
G_T(\Theta)$, by Theorem \ref{thm-adm-sim}, there exists a
$\xvt\in\Theta$ such that $f=(\xvt\bullet S)_T$ and $\xvt\bullet
S$ is a uniformly integrable $\Q$-martingale for each
$\Q\in\Mc^e$. We shall show $\xvt$ satisfies
\reff{cond-admissible-s}. The first line of
\reff{cond-admissible-s} is obvious.

By Theorem \ref{thm-adm-sim}, there exists an approximating sequence
$(\xvt^j)\subset\Theta^s$ for $\xvt$ such that $(\xvt^j\bullet
S)_T\to(\xvt\bullet S)_T$ in $L^2(\P)$, whence for any $g\in\Dc^s$,
$(\xvt^j\bullet S)_TZ^g_T\to(\xvt\bullet S)_TZ^g_T$ in $L^1(\P)$.
Moreover, Lemma \ref{lma-xth-s} below yields that, for any
$g\in\Dc^s$,
\begin{eqnarray}\label{lim-vt-1}
(\xvt^j\bullet S)_tZ^g_t=\E[(\xvt^j\bullet
S)_TZ^g_T|\Fc_t]\overset{L^1(\P)}{-\!-\!\!\!\longrightarrow}\E[(\xvt\bullet
S)_TZ^g_T|\Fc_t].\end{eqnarray} Then \reff{lim-vt-1} and
\reff{lim-vt-2} imply $\E[(\xvt\bullet S)_TZ^g_T|\Fc_t]=(\xvt\bullet
S)_tZ^g_t$ for each $g\in\Dc^s$ and therefore $\xvt$ satisfies the
second line of \reff{cond-admissible-s}. \qed

The following lemma has been used to prove Theorem
\ref{thm-theta-u}.

\begin{Lemma}\label{lma-xth-s} Under assumptions (H1) and (H2),
$\Theta^s\subset\Theta^u$.\end{Lemma}

{\bf Proof.} Let $\xvt\in\Theta^s$. The first line of
\reff{cond-admissible-s} is obvious. Let $\xvt\in\Theta^s$ have form
\reff{trade-simple}. For any $0\le s\le t\le T$ and $A\in\Fc_{s}$,
we can see
$$\widetilde{\xvt}:=(\id_A\id_{]\!]s,t]\!]})\xvt
=\sum_{i=1}^n\widetilde h_i\id_{]\!]\widetilde T_{1i},\widetilde
T_{2i}]\!]}\in\Theta^s,$$ where
$$\widetilde T_{1i}=(T_{1i}\vee s)\wedge(T_{2i}\wedge t),\quad
\widetilde T_{2i}=T_{2i}\wedge t,\quad
 \widetilde h_i=h_i\id_A\id_{[(T_{1i}\vee s)\le \widetilde T_{2i}]}.$$
Obviously,
$$(\widetilde\xvt\bullet S)_T=(\id_A\id_{]\!]s,t]\!]})\bullet(\xvt\bullet S)
=((\xvt\bullet S)_{t}-(\xvt\bullet S)_{s})\id_A.$$ For any
$g\in\Dc^s$, we have
\begin{eqnarray*}
\E[(\xvt\bullet S)_{t}Z^g_{t}\id_A]-\E[(\xvt\bullet
S)_{s}Z^g_{s}\id_A]
&=&\E[(\xvt\bullet S)_{t}Z^g_{T}\id_A]-\E(\xvt\bullet S)_{s}Z^g_{T}\id_A]\\
&=&\E[(\widetilde\xvt\bullet S)_Tg]\\
&=&0,\end{eqnarray*} which implies $\E[(\xvt\bullet
S)_{t}Z^g_{t}|\Fc_{s}]= (\xvt\bullet S)_{s}Z^g_{s}$ a.s. and
therefore $(\xvt\bullet S)Z^g$ is a uniformly integrable martingale.
That is just the second line of \reff{cond-admissible-s}.\qed

%

\section{Variance-optimal signed martingale measure}\label{sec-vsmm}
\setcounter{equation}{0} \setcounter{Assumption}{0}
\setcounter{Theorem}{0} \setcounter{Proposition}{0}
\setcounter{Corollary}{0} \setcounter{Lemma}{0}
\setcounter{Definition}{0} \setcounter{Remark}{0}

It is easy to see that, under hypothesis (H1) and that
$\Mc^s\ne\emptyset$, $\Dc^s$ is a non-empty closed convex subset of
$L^2(\P)$. Therefore there exists a unique $\Q^*\in\Mc^s$ such that
$g^*:={d\Q^*\over d\P}$ has minimal $L^2(\P)$-norm in $\Dc^s$. The
signed measure $\Q^*$ is called the variance-optimal signed
martingale measure (VSMM, for short) for $S$, since $g^*$ minimizes
$\var[g]=\E[g^2]-1$ over $g\in\Dc^s$.
By Lemma 2.1 of Delbaen and Schachermayer (1996a) and Lemma
\ref{lma-opt-1} below, we have:

\begin{Lemma}\label{lma-ds-c}
If $\Mc^s\ne\emptyset$, then $g^*$ is the
    unique element of $\overline{G_T(\Theta^s)}+\R$ such that (as a linear functional on
    $L^2(\P)$) vanishing on $\overline{G_T(\Theta^s)}$ and equaling $1$
    on the constant function $1$.
\end{Lemma}

\begin{Lemma}\label{lma-opt-1}
If $\Mc^s\ne\emptyset$, then
$\overline{G_T(\Theta^s)+\R}=\overline{G_T(\Theta^s)}+\R$.
\end{Lemma}

{\bf Proof.} The ``$\supset$" part is clear. Let
$f\in\overline{G_T(\Theta^s)+\R}$, then there exist sequences
$(f^j)\subset G_T(\Theta^s)$ and $(\xd^j)\subset\R$ such that
$f^j+\xd^j\to f$ in $L^2(\P)$. For any $g\in\Dc^s$, we have
$\xd^j=\E[(f^j+\xd^j)g]\to\E[fg]$ and therefore $f^j\to f-\E[fg]$ in
$L^2(\P)$, which yields $f\in\overline{G_T(\Theta^s)}+\R$. Thus
``$\subset$" part also holds.\qed

Under assumptions (H1) and (H2), by Lemma \ref{lma-ds-c} and
Theorems \ref{thm-adm-sim}--\ref{thm-theta-u}, there exist
$\vartheta^*\in\Theta^u$ and $a\in\R$ such that
$g^*=a+(\xvt^*\bullet S)_T$. Then we have
$$\E[(g^*)^2]=\E[(a+(\xvt^*\bullet S)_T)g^*]=a$$ and therefore
\begin{eqnarray}\label{eq-gopt-xvt}
g^*=\E[(g^*)^2]+(\xvt^*\bullet S)_T.\end{eqnarray} For any fixed
$\Q\in\Mc^e$, let $\widetilde Z^*$ be the RCLL version of the
$\Q$-martingale defined by
\begin{eqnarray}\label{def-zopt}\widetilde Z^*_t=\E_\Q[g^*|\Fc_t],\quad t\in[0,T],\end{eqnarray}
then by \reff{eq-gopt-xvt} and $\vartheta^*\in\Theta^u$, we have
\begin{eqnarray}\label{eq-zopt-xvt}\widetilde
Z^*_t=\E[(g^*)^2]+(\xvt^*\bullet S)_t,\quad t\in[0,T].\end{eqnarray}
Thus the definition of $\widetilde Z^*$ is independent of the choice
of $\Q\in\Mc^e$. Moreover, for each $g\in\Dc^s$, $\widetilde Z^*
Z^g$ is a uniformly integrable martingale since $\xvt^*\in\Theta^u$.
The above arguments lead to the following lemma, which extends Lemma
2.2 of Delbaen and Schachermayer (1996a).

\begin{Lemma}\label{lma-zopt} Under assumptions (H1) and (H2), we
have:
\begin{description}
    \item[(a)] $\widetilde Z^*$, as defined in \reff{def-zopt}, is
independent of the choice of $\Q\in\Mc^e$;
    \item[(b)] There exists $\vartheta^*:=(\vartheta^{*1},\dots,\vartheta^{*d})^\tr\in\Theta^u$ such that \reff{eq-zopt-xvt} holds;
    \item[(c)] For each $g\in\Dc^s$, $\widetilde Z^* Z^g$ is a
uniformly integrable martingale. In particular, $\widetilde Z^* Z^*$
is a uniformly integrable martingale, where $Z^*$ is the RCLL
version of the martingale defined by $Z^*_t=\E[g^*|\Fc_t]$.
\end{description}
\end{Lemma}

\begin{Lemma}\label{lma-zz-ne0}
Under assumptions (H1) and (H2), we have:\begin{description}
    \item[(a)] $\widetilde Z^* Z^*\ge0$;
    \item[(b)] Let $\tau$ be the first time when $\widetilde Z^*Z^*$ hits $0$, that is,
$$\tau=\inf\{t\in[0,T]: \widetilde Z^*_tZ^*_t=0\}\mbox{ with
}\inf\emptyset=\infty,$$ then  $\widetilde
Z^*Z^*\id_{[\![\tau,T]\!]}=0$;
    \item[(c)] If $\P(g^*=0)=0$, then $\widetilde Z^* Z^*>0$.
\end{description}
\end{Lemma}

{\bf Proof.} By Lemma \ref{lma-zopt}(c), we have for any
$t\in[0,T]$ that $\widetilde Z^*_t Z^*_t=\E[(g^*)^2|\Fc_t]\ge0$
a.s., and therefore (a) holds since $\widetilde Z^* Z^*$ is RCLL.
Obviously, (a) and (b) imply (c). It remains to show (b).

For any stopping time $\xs$ with $0\le\xs\le T$ a.s., Doob
stopping theorem leads to
$$\E[\widetilde Z^*_\xs Z^*_\xs]=\E[\widetilde Z^*_{\xs\land \tau}Z^*_{\xs\land \tau}]
=\E[\widetilde Z^*_\xs Z^*_\xs\id_{[\xs<\tau]}]$$ and hence
$\E[X_\xs]=\E[\widetilde Z^*_\xs Z^*_\xs\id_{[\xs\ge\tau]}]=0$,
where $X=\widetilde Z^*Z^*\id_{[\![\tau,T]\!]}$. Then by section
theorem, $X=0$ indistinguishably, that is just (b). \qed

%
%

\section{Mean-variance hedging}\label{sec-mv}
\setcounter{equation}{0} \setcounter{Assumption}{0}
\setcounter{Theorem}{0} \setcounter{Proposition}{0}
\setcounter{Corollary}{0} \setcounter{Lemma}{0}
\setcounter{Definition}{0} \setcounter{Remark}{0}

In this section, we always assume (H1) and (H2). The mean-variance
hedging problem is, for any $H\in L^2(\P)$,  to
\begin{eqnarray}\label{prob-mvth}
\mbox{Minimize }\E[(H-G_T(\xvt))^2]\quad \mbox{ subject to
}\xvt\in\xT.\end{eqnarray}

\subsection{$L^2(\P)$-orthogonal decomposition}

Recalling Theorem \ref{thm-adm-sim}(b) and Lemma \ref{lma-opt-1}, we
know $G_T(\Theta)=\overline{G_T(\Theta^s)}$ and
$G_T(\Theta)+\R=\overline{G_T(\Theta^s)+\R}$. So both $G_T(\Theta)$
and $G_T(\Theta)+\R$ are closed subspaces of $L^2(\P)$. Thus
\reff{prob-mvth} always allows a solution. By Lemma
\ref{lma-ds-1}(b), we know $$\Dc^s-g^*:=\{g-g^*: g\in\Dc^s\}$$ is a
closed subspace of $L^2(\P)$ and by Lemma \ref{lma-bipolar},
$G_T(\Theta)+\R=(\Dc^s-g^*)^\perp$. Hereafter, $\{\dots\}^\perp$
stands for the orthogonal complement in $L^2(\P)$. We denote by
$\pi$ the projection in $L^2(\P)$ on $G_T(\Theta)^\perp$, then
$\pi(1)\in G_T(\Theta)+\R$ and by Lemma \ref{lma-ds-c},
$g^*={\pi(1)\over\E[\pi(1)]}$. The above arguments lead to an
orthogonal decomposition of $G_T(\Theta)+\R$ as follows:
$$G_T(\Theta)+\R=G_T(\Theta)\oplus g^*\R,$$
where $g^*\R$ is the linear space spanned by $g^*$, that is,
$g^*\R=\{\xa g^*: \xa\in\R\}$. Moreover, $L^2(\P)$ can be
orthogonally decomposed as follows:
$$L^2(\P)=G_T(\Theta)\oplus g^*\R\oplus(\Dc^s-g^*).$$
Consequently, by Theorem \ref{thm-theta-u}, we have the following
easy theorem.

\begin{Theorem}\label{thm-easy} Assume (H1) and (H2).
Any $H\in L^2(\P)$ admits a unique orthogonal decomposition
\begin{eqnarray}\label{dcp-h}
H=G_T(\vartheta^H)\oplus\xa^Hg^*\oplus (g^H-g^*),\end{eqnarray}
where
$\vartheta^H:=(\vartheta^{H1},\dots,\vartheta^{Hd})^\tr\in\Theta^u$,
$\xa^H={\E[Hg^*]\over\E[(g^*)^2]}$ and $g^H\in\Dc^s$. Moreover,
$\vartheta^H$ solves \reff{prob-mvth}.\end{Theorem}

\subsection{Num\'eraire approach}

Hereafter, in addition to (H1) and (H2), we always assume
\begin{description}
    \item[(H3)] $\P(g^*=0)=0$.
\end{description}
By Lemma \ref{lma-zz-ne0}, $\wt Z^*Z^*>0$.

Following GLP 1998, see also RS 1997 and Arai (2005), we define an
$\R^{d+1}$-valued process $Y$ and a new probability measure
$\widetilde\P$ as follows:
\begin{eqnarray*}
Y^0&:=&(\widetilde Z^*)^{-1},\\
Y^i&:=&S^i(\widetilde Z^*)^{-1},\quad i=1,\dots,d,\\
{d\widetilde\P\over d\P}&:=&{(g^*)^2\over\E[(g^*)^2]}={(\widetilde
Z^*_T)^2\over\E[(g^*)^2]}.
\end{eqnarray*}
It is clear that, for any $H\in L^2(\P)$ and $\xvt\in\Theta$,
\begin{eqnarray}\label{norm-p-p}
\E[(H-G_T(\xvt))^2]=\E[(g^*)^2]\cdot\E_{\widetilde\P}\left[\left({H\over\widetilde
Z^*_T}-{G_T(\xvt)\over\widetilde
Z^*_T}\right)^2\right].\end{eqnarray}

The following notations will be used:
\begin{itemize}
    \item The space $\Mcc(\widetilde\P)$ (resp. $\Mcc_{\mathrm{loc}}(\widetilde\P)$)
    consists of all uniformly integrable (resp. local) $\widetilde\P$-martingales;
    \item The space $\Mcc^2(\widetilde\P)$ (resp. $\Mcc^2_{\mathrm{loc}}(\widetilde\P)$)
    consists of all (resp. locally) square-integrable $\widetilde\P$-martingales and
$\Mcc^2_0(\widetilde\P)=\{M\in\Mcc^2(\widetilde\P): M_0=0\}$;
    \item The space $\Psi$ consists of all processes $\psi\in\Lc(Y)$ such that $\psi\bullet
Y\in\Mcc^2(\widetilde\P)$;
    \item  The space $\wt\Theta$ consists of all processes $\vartheta\in\Lc(S)$ satisfying
$\int_0^T\xvt_t\,dS_t\in L^2(\P)$ and $(\xvt\bullet S)Z^*$ is a
uniformly integrable martingale.
\end{itemize}
Lemma \ref{lma-xth-s} shows
$\Theta^s\subset\Theta^u\subset\wt\Theta$.

\begin{Proposition}\label{prop-s-y}
Under assumptions (H1), (H2) and (H3), we have
\begin{eqnarray}\label{gt-y}{1\over\widetilde Z^*_T}\,G_T(\Theta^u)=\{(\psi\bullet Y)_T:
\psi\in\Psi\}.\end{eqnarray} Moreover, the relation between
$\xvt\in\Theta^u$ and $\psi\in\Psi$ is given by
\begin{eqnarray*}
\psi^i&=&\xvt^i,\quad i=1,\dots,d,\\
\psi^0&=&(\xvt\bullet S)-\xvt^\tr S
\end{eqnarray*} and
\begin{eqnarray}\label{eq-xvt-psi}
\xvt^i=\psi^i+\xvt^{*i}(\psi\bullet Y-\psi^\tr Y),\quad
i=1,\dots,d,\end{eqnarray} where $\vartheta^*$ is defined as in
Lemma \ref{lma-zopt}(b).\end{Proposition}

{\bf Proof.} The relation between $\xvt$ and $\psi$ and the fact
\begin{eqnarray}\label{xvt-y}
    \{\xvt\bullet S: \xvt\in\Lc(S)\}
=\{(\psi\bullet Y)\widetilde Z^*: \psi\in\Lc(Y)\}
\end{eqnarray}
has been proved in Proposition 8 of RS 1997.

In order to prove the ``$\subset$" part of \reff{gt-y}, it is enough
to show the corresponding $\psi$ is in $\Psi$ if
$\vartheta\in\wt\Theta$. Actually, let $\vartheta\in\wt\Theta$, then
$(\vartheta\bullet S)_T\in L^2(\P)$ and $(\xvt\bullet S)Z^*$ is a
uniformly integrable martingale.  By \reff{xvt-y}, $\xvt\bullet
S=(\psi\bullet Y)\wt Z^*$ and therefore $(\psi\bullet Y)\widetilde
Z^*Z^*$ is a uniformly integrable martingale. On the other hand,
Lemma \ref{lma-zopt}(c) shows
\begin{eqnarray}\label{dpp-dp}\E\left[\left.{d\widetilde\P\over
d\P}\right|\Fc_t\right]={\widetilde
Z^*_tZ^*_t\over\E[(g^*)^2]},\end{eqnarray} and therefore
$(\psi\bullet Y)\in\Mcc(\wt\P)$. Moreover,
$$\E_{\widetilde\P}[(\psi\bullet Y)_T^2]={\E[(\xvt\bullet S)_T^2]\over \E[(g^*)^2]}<\infty.$$
Thus $\psi\in\Psi$, which implies the ``$\subset$" part of
\reff{gt-y}.

It is worth noting that the argument in the previous paragraph
yields particularly that $Y^i\in\Mcc^2_{\mathrm{loc}}(\widetilde
\P)$ for each $i=0,1,\dots,d$, since $S\in\Lc^2_{\mathrm{loc}}(\P)$
and $\Theta^s\subset\wt\Theta$. On the other hand, by Lemma
\ref{lma-zopt}(c) , for each $g\in\Dc^s$, $Z^g\wt Z^*$ is a
uniformly integrable martingale, and then by \reff{dpp-dp},
$Z^g(Z^*)^{-1}\in\Mcc(\wt\P)$. Obviously,
$$\E_{\wt\P}[(Z^g_T(Z^*_T)^{-1})^2]={\E[g^2]\over\E[(g^*)^2]}<\infty,$$
thus $Z^g(Z^*)^{-1}\in\Mcc^2(\wt\P)$. Moreover, by
$\Theta^s\subset\Theta^u$ and by $S\in\Lc^2_{\mathrm{loc}}(\P)$, we
can see, for each $i=1,\dots,d$ and each $g\in \Dc^s$, $S^iZ^g$ is a
local martingale, that is,
$Y^iZ^g(Z^*)^{-1}\in\Mcc_{\mathrm{loc}}(\wt\P)$. Similarly, for each
$g\in\Dc^s$, $Z^g$ is a square-integrable martingale, that is,
$Y^0Z^g(Z^*)^{-1}\in\Mcc^2(\wt\P)$. To conclude this paragraph, we
have $\langle Y^i,Z^g(Z^*)^{-1}\rangle(\wt\P)=0$, for each
$i=0,1,\dots,d$ and each $g\in\Dc^s$.

Now we are in a position to prove the ``$\supset$" part. Let
$\psi\in\Psi$ and the corresponding $\xvt$ be given by
\reff{eq-xvt-psi}, we should show $\vartheta\in\Theta^u$. Actually,
by the results in the previous paragraph, we have for each
$g\in\Dc^s$ that $\langle\psi\bullet
Y,Z^g(Z^*)^{-1}\rangle(\wt\P)=0$ and therefore $(\psi\bullet
Y)Z^g(Z^*)^{-1}\in\Mcc(\wt\P)$. That is, by \reff{xvt-y} ,
$(\vartheta\bullet S)(\wt Z^*)^{-1}Z^g(Z^*)^{-1}\in\Mcc(\wt\P)$, and
then by \reff{dpp-dp}, $(\vartheta\bullet S)Z^g$ is a uniformly
integrable martingale, for each $g\in\Dc^s$. Now we have shown
$\vartheta\in\Theta^u$.\qed

Actually, in the proof of the previous proposition, we have
implicitly shown the following
\begin{Corollary}\label{cor-xvt-p}
Under assumptions (H1), (H2) and (H3), we have $\wt\Theta=\Theta^u$.
\end{Corollary}

\begin{Remark}\label{rmk-xvt-p}{\rm
For continuous semimartingale models, $\Theta^u$ was defined as the
working space of admissible strategies in GLP 1998 and $\wt\Theta$
was used by RS 1997. In the continuous case, since $\Q^*\in\Mc^e$,
the previous corollary shows $\Theta^u=\wt\Theta$. For the
discontinuous case, when $\Q^*\in \Mc^e$, $\wt\Theta$ was used in
Arai (2005). Anyway, we have $\Theta^u=\wt\Theta$ if (H3) is
satisfied. }\end{Remark}

In view of Proposition \ref{prop-s-y} and \reff{norm-p-p},
\reff{prob-mvth} is equivalent to the problem to
\begin{eqnarray}\label{min-pp}
\mbox{Minimize } \E_{\widetilde\P}\left[\left({H\over\widetilde
Z^*_T}-\int_0^T \psi\, dY\right)^2\right]\quad \mbox{over }
\psi\in\Psi.\end{eqnarray} The solution of \reff{min-pp} is given by
the GKW decomposition
\begin{eqnarray}\label{dcp-gwk}
{H\over\widetilde Z^*_T}=\E_{\widetilde\P}\left[{H\over\widetilde
Z^*_T}\right]+\int_0^T\psi^H_t\, dY_t+L^H_T,\end{eqnarray} under
$\widetilde\P$, where
$\psi^H:=(\psi^{H0},\psi^{H1},\dots,\psi^{Hd})^\tr\in\Psi$ and
$L^H\in\Mcc^2_0(\widetilde\P)$ is strongly
$\widetilde\P$-orthogonal to $Y$.

The following theorem, which extends the corresponding results of
GLP 1998, RS 1997 and Arai (2005) to our situation, shows the
relation between decompositions \reff{dcp-h} and \reff{dcp-gwk}.

\begin{Theorem}\label{thm-relation-1} Assume (H1), (H2) and (H3).
For decompositions \reff{dcp-h} and \reff{dcp-gwk}, we have
\begin{eqnarray}\label{eq-relation-1-1}
\xvt^{Hi}&=&\psi^{Hi}+\xvt^{*i}(\psi^H\bullet Y-(\psi^H)^\tr
Y),\quad i=1,\dots,d,\\
g^H-g^*&=&L^H_Tg^*.\label{eq-relation-1-2}
\end{eqnarray}
\end{Theorem}

{\bf Proof.} Obviously, $\E_{\widetilde\P}\left[{H\over\widetilde
Z^*_T}\right]={\E[Hg^*]\over\E[(g^*)^2]}=\xa^H$. Proposition
\ref{prop-s-y} and \reff{norm-p-p} imply \reff{eq-relation-1-1},
since $\psi^H$ solves \reff{min-pp}. By \reff{dcp-gwk}, we have
\begin{eqnarray*}
H&=&\widetilde Z^*_T\int_0^T\psi^H_t\,
dY_t+\E_{\widetilde\P}\left[{H\over\widetilde
Z^*_T}\right]\widetilde Z^*_T+L^H_T\widetilde Z^*_T\\
&=&G_T(\vartheta^H)+\xa^Hg^*+L^H_Tg^*,\end{eqnarray*} which, by
\reff{dcp-h}, implies \reff{eq-relation-1-2}.\qed

\subsection{Decompositions under VSMM}

In this subsection we study two kinds of decompositions under the
VSMM $\Q^*$: one is directly derived from decompositions
\reff{dcp-h} and the other one is the GKW decomposition. Under
assumption (H3), by Lemma \ref{lma-zz-ne0}, $\wt Z^*Z^*>0$
indistinguishably. For any $H\in L^2(\P)$, we can define $V^H$ as
follows:
$$V^H_t={\E[HZ^*_T|\Fc_t]\over Z^*_t}\quad\mbox{ for all
}t\in[0,T].$$ Obviously, $V^HZ^*$ is a uniformly integrable
martingale. If $\Q^*\in\Mc^e$, then $V^H_t=\E_{\Q^*}[H|\Fc_t]$.

In the following theorem, we introduce a decomposition under $\Q^*$
which is directly derived from decompositions \reff{dcp-h}.
\begin{Theorem}\label{thm-relation-2} Assume (H1), (H2) and (H3).
For any $H\in L^2(\P)$, $V^H$ has a unique decomposition
\begin{eqnarray}\label{dcp-vh}
V^H=V^H_0+\varphi^H\bullet S+K^H,\end{eqnarray} where
$V^H_0=\E[Hg^*]$, $\varphi^H\in\Theta^u$, $K^H_0=0$, $K^HZ^*$ is a
uniformly integrable martingale and $K^H_T\in
(G_T(\Theta)+\R)^\perp$. The relation between decompositions
\reff{dcp-h} and \reff{dcp-vh} is:
\begin{eqnarray}\label{dcp-relation-2-1}
\varphi^H&=&\vartheta^H+\xa^H\vartheta^*,\\
K^H_t&=&{\E[(g^H-g^*)Z^*_T|\Fc_t]\over
Z^*_t}.\label{dcp-relation-2-2}
\end{eqnarray}
Moreover, we have
\begin{eqnarray}\label{dcp-relation-2-3}
K^H=L^H\wt Z^*,\end{eqnarray} where $L^H$ is given by
decomposition \reff{dcp-gwk}.
\end{Theorem}

{\bf Proof.} ``Uniqueness": In addition to \reff{dcp-vh}, suppose
$V^H$ has another decomposition $V^H=V^H_0+\varphi\bullet S+K$ with
$\varphi$ and $K$ satisfying the corresponding conditions of
$\varphi^H$ and $K^H$ respectively, then $(\varphi-\varphi^H)\bullet
S=K^H-K$ and therefore $K^H_T-K_T\in G_T(\Theta)$. On the other
hand, the decompositions require $K^H_T-K_T\in
(G_T(\Theta)+\R)^\perp$. Thus $K^H_T=K^T$ a.s. Since both $K^HZ^*$
and $KZ^*$ are uniformly integrable martingales, we have further
that $K^H=K$ and hence $\varphi^H\bullet S=\varphi\bullet S$. This
completes the uniqueness.

``Existence": For any $H\in L^2(\P)$, $H$ has decomposition
\reff{dcp-h}. By $\vartheta^H\Theta^u$ and Lemma \ref{lma-zopt}(c),
we have
$$\E[HZ^*_T|\Fc_t]=G_t(\vartheta^H)Z^*_t+\xa^H\wt
Z^*_tZ^*_t+\E[(g^H-g^*)g^*|\Fc_t],$$ that is, by
\reff{eq-zopt-xvt},
$$V^H_t=G_t(\vartheta^H)+\xa^H\wt
Z^*_t+K^H_t=\E[Hg^*]+G_t(\varphi^H)+K^H_t,$$ where $\varphi^H$ and
$K^H$ are respectively given by \reff{dcp-relation-2-1} and
\reff{dcp-relation-2-2}. It is easy to verify that
$\varphi^H\in\Theta^u$, $K^H_0=0$, $K^HZ^*$ is a uniformly
integrable martingale and $H^H_T=g^H-g^*\in(G_T(\Theta)+\R)^\perp$.

Finally, by \reff{dpp-dp}, $L^H\widetilde Z^*Z^*$ is a uniformly
integrable martingale since $L^H\in\Mcc^2_0(\widetilde\P)$. Then
\reff{eq-relation-1-2} and \reff{dcp-relation-2-2} yield
\reff{dcp-relation-2-3}. \qed

\begin{Remark}{\rm
In decomposition \reff{dcp-vh}, since
$K^H_T\in(G_T(\Theta)+\R)^\perp$, $(V^H_0,\varphi^H)$ solves the
problem to
\begin{eqnarray*}
\mbox{Minimize }\E[(H-x-(\varphi\bullet S)_T)^2]\quad\mbox{over
}(x,\varphi)\in\R\times\Theta.\end{eqnarray*}
 }\end{Remark}

Now we consider the GKW decompositions of local martingales under
$\Q^*$. To this end, we assume further that $\Q^*\in\Mc^e$.

The GKW decomposition of $V^H$ is:
\begin{eqnarray}\label{dcp-gkw-vh}
V^H=V^H_0+\eta^H\bullet S+N^H,\end{eqnarray} where $V^H_0=\E[Hg^*]$,
$\eta^H\in\Lc(S)$, and $N^H$ satisfies $N^H_0=0$ and both $N^H$ and
$N^HS$ are local $\Q^*$-martingales.

On the other hand, by Theorem \ref{thm-relation-2}, $V^H$ can be
decomposed as follows:
$$V^H=V^H_0+(\vartheta^H+\xa^H\vartheta^*)\bullet S+L^H\wt Z^*.$$
Let \begin{eqnarray}\label{eq-def-llh}J^H:=L^H\wt
Z^*-L^H_-\bullet\wt Z^*,\end{eqnarray} then by \reff{eq-zopt-xvt},
we have
\begin{eqnarray}V^H=V^H_0+(\vartheta^H+(\xa^H+L^H_-)\vartheta^*)\bullet
S+J^H,\label{eq-vt-0}
\end{eqnarray}
Obviously, $J^H$ is a local $\Q^*$-martingale with $J^H_0=0$.

Assume the GKW decomposition of $J^H$ is
\begin{eqnarray}\label{dcp-gkw-j}
J^H=\eta^J\bullet S+N^J,\end{eqnarray} where $\eta^J\in\Lc(S)$,
$N^J_0=0$, and both $N^J$ and $N^JS$ are local $\Q^*$-martingales.
In view of \reff{dcp-gkw-vh}, \reff{eq-vt-0} and \reff{dcp-gkw-j},
by the uniqueness of the GKW decomposition, we have $N^H=N^J$ and
\begin{eqnarray}\label{eq-eta-hj}
\eta^H=\vartheta^H+(\xa^H+L^H_-)\vartheta^*+\eta^J.\end{eqnarray}
The above arguments lead to the following theorem, which extends
Theorem 4.1 of Arai (2005) to our settings.

\begin{Theorem}\label{thm-feedback}
Assume (H1), (H2) and $\Q^*\in\Mc^e$. Let $H\in L^2(\P)$, if $V^H$
and $J^H$ allow the GKW decompositions \reff{dcp-gkw-vh} and
\reff{dcp-gkw-j} respectively, then the solution $\vartheta^H$ of
\reff{prob-mvth} satisfies the following feedback equation:
\begin{eqnarray}\label{eq-feedback}
\vartheta^H=\eta^H-\eta^J-{\vartheta^*\over\wt
Z^*_-}\,(V^H_--G_-(\vartheta^H)).\end{eqnarray}
\end{Theorem}

{\bf Proof.} By integration by parts and \reff{eq-zopt-xvt}, we have
\begin{eqnarray*}
\wt Z^*(\xa^H+L^H)&=&\xa^H\wt
Z^*_0+((\xa^H+L^H_-)\vartheta^*)\bullet S+\wt Z^*_-\bullet
L^H+[\wt Z^*,L^H]\\
&=&V^H-G(\vartheta^H)\qquad\qquad(\mbox{by
\reff{eq-vt-0}})\end{eqnarray*} and therefore
$$\xa^H+L^H_-={V^H_--G_-(\vartheta^H)\over\wt Z^*_-}.$$
Then by \reff{eq-eta-hj},
$$\vartheta^H=\eta^H-\eta^J-(\xa^H+L^H_-)\vartheta^*=\eta^H-\eta^J-{\vartheta^*\over\wt
Z^*_-}\,(V^H_--G_-(\vartheta^H)).$$\qed

\begin{Proposition}\label{prop-etaj-0} Under the conditions of
Theorem \ref{thm-feedback}, $\eta^J=0$ if and only if $\sum \xD
L^H\xD\wt Z^*\xD S$ is a local $\Q^*$-martingale. If it is the case,
then the solution $\vartheta^H$ of \reff{prob-mvth} satisfies the
following feedback equation:
\begin{eqnarray}\label{eq-feedback-0}
\vartheta^H=\eta^H-{\vartheta^*\over\wt
Z^*_-}\,(V^H_--G_-(\vartheta^H)),\end{eqnarray} where $\eta^H$ is
given by the GKW decomposition \reff{dcp-gkw-vh}.
\end{Proposition}

{\bf Proof.} Obviously, in decomposition \reff{dcp-gkw-j},
$\eta^J=0$ if and only if $[J^H,S]$ is a local $\Q^*$-martingale. By
integration by parts, one can compute that
\begin{eqnarray*}
[J^H,S]=[\wt Z^*_-\bullet L^H,S]+[[L^H,\wt Z^*],S]=\wt
Z^*_-\bullet[L^H,S]+\tsm\,\xD L^H\xD\wt Z^*\xD S.\end{eqnarray*} By
Lemma \ref{lma-lh-s} below, $\wt Z^*_-\bullet[L^H,S]$ is a local
$\Q^*$-martingale, and therefore 
the conclusion of the proposition follows.\qed

\begin{Remark}{\rm If $S$ is a continuous semimartingale satisfying (H1) and (H2),
then $\Q^*\in\Mc^e$ automatically holds. In this case, due to the
continuity of $S$, it always holds that $\sum \xD L^H\xD\wt Z^*\xD
S=0$ and therefore Proposition 10 of RS 1997 can be recovered here.}
\end{Remark}

\begin{Remark}{\rm In general, however, it is
rather restrictive to assume $\sum \xD L^H\xD\wt Z^*\xD S$ is a
local $\Q^*$-martingale. This fact was observed by Arai
(2005).}\end{Remark}

The following lemma has been used to prove Proposition
\ref{prop-etaj-0}.

\begin{Lemma}\label{lma-lh-s}
Under the conditions of Theorem \ref{thm-feedback}, for any $H$,
$[L^H,S]$ is a local $\Q^*$-martingale.
\end{Lemma}

{\bf Proof.} Since $L^H$ is strongly $\wt\P$-orthogonal to $Y$, we
know $L^HY^i\in\Mcc(\wt\P)$ for each $i=0,\dots,d$. For $i=0$,
$L^HY^0=L^H(\wt Z^*)^{-1}$ and therefore $L^H$ is a uniformly
integrable $\Q^*$-martingale. For each $i=1,\dots,d$,
$L^HY^i=L^HS^i(\wt Z^*)^{-1}$ and therefore $L^HS^i$ is a uniformly
integrable $\Q^*$-martingale. On the other hand, $S$ is a
$\Q^*$-martingale and therefore we can see $[L^H,S]$ is a local
$\Q^*$-martingale.\qed

%
%
%
%
%
%
%
%

%
%

\newpage

\vspace{10mm} \noindent{\Large\bf References} {
\begin{description}

\item Arai, T. (2005): An extension of mean-variance hedging to
the discontinuous case, {\it Finance and Stochastics} {\bf 9},
129-139.

\item \v{C}ern\'y, A., J. Kallsen (2005): On the structure of
general mean-variance hedging strategies, preprint.

\item Choulli, T., L. Krawczyk and C. Stricker (1998): $\cal
E$-martingales and their applications in mathematical finance,
{\it Annals of Probability} {\bf 26}, 853-876.

\item Choulli, T., C. Stricker and L. Krawczyk (1999): On
Fefferman and Burkholder-Davis-Gundy inequalities for $\cal
E$-martingales, {\it Probability Theory and Related Fields} {\bf
113}, 571-597.

\item Delbaen, F., P. Monat, W. Schachermayer, M. Schweizer and C.
Stricker (1997): Weighted norm inequalities and hedging in
incomplete markets, {\it Finance and Stochastics} {\bf 1},
181-227.

\item Delbaen, F. and W. Schachermayer (1996a): The
Variance-Optimal Martingale Measure for Continuous Processes, {\it
Bernoulli} {\bf 2}, 81-105.

\item Delbaen, F. and W. Schachermayer (1996b): Attainable claims
with $p$'th moments, {\it Ann. Inst. H. Poincar\'e Probab.
Statist.} {\bf 32}, 743--763.


\item Gouri\'eroux, C., J. P. Laurent and H. Pham (1998):
Mean-Variance Hedging and Num\'eraire, {\it Mathematical Finance}
{\bf 8}, 179-200.

\item Grandits, P. and L. Krawczyk (1998): Closedness of some
spaces of stochastic integrals, {\it S\'eminaire de Probabilit\'es
XXXII}, LN in Math {\bf 1686}, 73-85. Berlin: Springer.

\item Hou, C. and I. Karatzas (2004): Least-squares approximation
of random variables by stochastic integrals, {\it Adv. Stud. Pure
Math.} {\bf 41}, 141-161.

\item Jacod, J. (1979): {\it Calcul stochastique et probl\'emes de
martingales}, LN in Math {\bf 714}. Berlin: Springer.

\item Jacod, J. (1980): Int\'egrles stochastiques par rapport
\`a une semimartingale vectorielle et changement de filtrations,
{\it S\'em. Probab. XIV, LN in Math.} {\bf 784}, 161-172, Springer.

\item Jacod, J. and A. Shiryaev (1987): {\it Limit Theorems for
Stochastic Process}, Berlin: Springer.





\item Rheinl\"ander, T. and M. Schweizer (1997): On
$L^2$-Projections on a Space of Stochastic Integrals, {\it Annals
of Probability} {\bf 25}, 1810-1831.

\item Schachermayer, W. (2005): A note on arbitrage and closed
convex cones, {\it Mathematical Finance} {\bf 15}, 183-189.



\item Schweizer, M. (1996): Approximation pricing and the
variance-optimal martingale measure, {\it Annals of Probability}
{\bf 24}, 206-236.


\item Shiryaev, A.N.  and A.S. Cherny (2002): Vector stochastic integrals and the
fundamental theorems of asset pricing, {\it Proceedings of the
Steklov Mathematical Institute} {\bf 237}, 12-56.


\item Xia, J. and J.A. Yan (2006): Markowitz's portfolio
optimization in an incomplete market, {\it Mathematical Finance}
{\bf 16}, 203-216.

%

\end{description}
}

\end{document}